Article: CJB/2010/11

# Eight circles through the Orthocentre

**Christopher J Bradley**

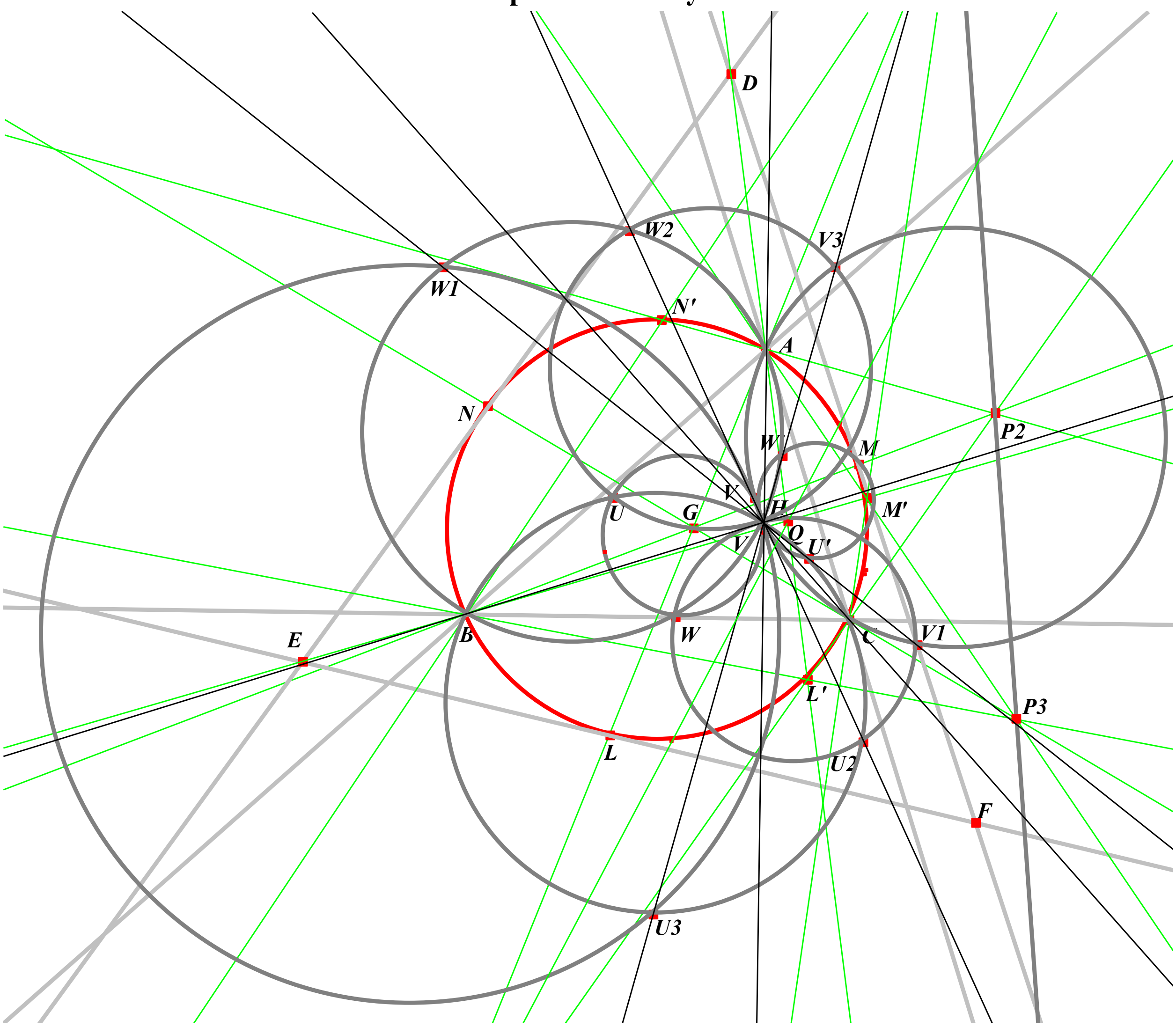

**Fig. 1**

**Circles BHC, CHA, AHB and five Hagge Circles**

## 1. Introduction

An appreciation of Fig. 1 (drawn with *CABRI II* plus) requires a description of how it was drawn and the order in which the various constructions were made. Start with a triangle ABC and its centroid G. AG, BG, CG meet the circumcircle Γ of ABC at points L, M, N respectively. Tangents at L, M, N to Γ are then drawn with those at M and N meeting at D. E and F are similarly defined. Γ then becomes an incircle of triangle DEF. Next DA, EB, FC are drawn and these turn out to meet at a point Q and when extended meet Γ again at L', M', N' respectively. Lines AL, BN', CM' now prove to be concurrent at $P_1$; BM, CL', AN' prove to be concurrent at $P_2$; and CN, AM', BL' prove to be concurrent at $P_3$. Furthermore it transpires that $P_1$, $P_2$, $P_3$ are collinear lying on the polar of Q with respect to Γ. Twelve more points are now obtained as reflections of L, M, N, L', M', N' in the sides of ABC. These are U, V, W the reflections of L, M, N respectively in the sides BC, CA, AB; U', V', W' are the reflections of L', M', N' respectively in the sides BC, CA, AB; $U_2$, $U_3$ are the reflections of M', N' respectively in the side BC; $V_3$, $V_1$ are the reflections of N', L' respectively in the side CA; and finally $W_1$, $W_2$ are the reflections of L', M' respectively in the side AB. Now five circles can be drawn, all of which pass through the orthocentre H, UVW is the Hagge circle of G, U'V'W' is the Hagge circle of Q, $UV_3W_2$ is the Hagge circle of $P_1$, $U_3VW_1$ is the Hagge circle of $P_2$ and $U_2V_1W$ is the Hagge circle of $P_3$. For the original account of Hagge circles see Hagge [1], for a modern account see Bradley and Smith [2]. Three more circles may now be drawn, circles BHC, CHA, AHB, and these circles remarkably contain the points U, U', $U_2$, $U_3$ and V, V', $V_3$, $V_1$ and W, W', $W_1$, $W_2$ respectively. It is also the case that lines $V_1W_1$, $W_2U_2$ and $U_3V_3$ all pass through H.

The analysis that follows in later sections proves all the above results using areal co-ordinates with ABC as triangle of reference. For an account of these co-ordinates and how to use them, see Bradley [3]. It must now be confessed that the point G was used in the construction and in the analysis that follows in order to provide an analysis that is as easy as possible, the areal co-ordinates of G being simply (1, 1, 1). The truly remarkable thing is that all the concurrences, collinearities, and dispositions of points on circles hold not just for the starting point G, but for a **any internal point of** ABC taken as starting point (except H itself and except the symmedian point K when in both cases the figure degenerates). The problem, of course, is that if one starts with P(l, m, n) the analysis becomes technically very difficult (because of the three extra variables) even with a good algebra computer package, but not only that, it becomes virtually unprintable. No excuses are therefore made for having taken the easiest option.

## 2. The points L, M, N, D, E, F

The equation of the line AG is y = z and this meets the circumcircle Γ with equation

$$a^2yz + b^2zx + c^2xy = 0 \tag{2.1}$$

at the point L with co-ordinates $L(-a^2, b^2 + c^2, b^2 + c^2)$. M and N follow by cyclic change and have co-ordinates $M(c^2 + a^2, -b^2, c^2 + a^2)$ and $N(a^2 + b^2, a^2 + b^2, -c^2)$.

The equation of the tangent to Γ at any point (d, e, f) is

$$a^2fy + a^2ez + b^2dz + b^2fx + c^2ex + c^2dy = 0. \tag{2.2}$$

It follows that the equation of the tangent at L is

$$(b^2 + c^2)^2x + a^2b^2y + c^2a^2z = 0. \tag{2.3}$$

The tangents at M and N are respectively

$$a^2b^2x + (c^2 + a^2)^2y + b^2c^2z = 0 \qquad (2.4)$$

and

$$c^2a^2x + b^2c^2y + (a^2 + b^2)^2z = 0. \qquad (2.5)$$

The tangents at M and N meet where both equations (2.4) and (2.5) hold. This is at the point D with co-ordinates D($-$ ($a^4 + a^2(b^2 + c^2) + 2b^2c^2$), $b^2(a^2 + b^2 - c^2)$, $c^2(c^2 + a^2 - b^2)$). The co-ordinates of E and F follow by cyclic change and are E($a^2(a^2 + b^2 - c^2)$, $-$ ($b^4 + b^2(c^2 + a^2) + 2c^2a^2$, $c^2(b^2 + c^2 - a^2)$) and F($a^2(c^2 + a^2 - b^2)$, $b^2(b^2 + c^2 - a^2)$, $-$ ($c^4 + c^2(a^2 + b^2) + 2a^2b^2$)).

## 3. The points Q, L', M', N', $P_1$, $P_2$, $P_3$

The critical result that makes the configuration possess so many interesting details is the concurrence of the lines DA, EB, FC at the point Q. This is the point that links L, M, N with the points L', M', N'. Without that link the figure would not contain enough circles to be of significance.

The equations of the lines DA, EB, FC are

$$c^2(c^2 + a^2 - b^2)y = b^2(a^2 + b^2 - c^2)z, \qquad (3.1)$$

$$a^2(a^2 + b^2 - c^2)z = c^2(b^2 + c^2 - a^2)x, \qquad (3.2)$$

$$b^2(b^2 + c^2 - a^2)x = a^2(c^2 + a^2 - b^2)y, \qquad (3.3)$$

respectively. These lines are concurrent at the point Q with co-ordinates (x, y, z), where $x = \frac{a^2}{b^2 + c^2 - a^2}, y = \frac{b^2}{c^2 + a^2 - b^2}, z = \frac{c^2}{a^2 + b^2 - c^2}$. Note that the isogonal of Q is the isotomic of H. Now the line DA with Equation (3.1), when extended, meets Γ, with Equation (2.1), at the point L', which therefore has co-ordinates L'($-$ 1, $2b^2/(c^2 + a^2 - b^2)$, $2c^2/(a^2 + b^2 - c^2)$). The co-ordinates of M' and N' follow by cyclic change and are M'($2a^2/(b^2 + c^2 - a^2)$, $-$ 1, $2c^2/(a^2 + b^2 - c^2)$) and N'($2a^2/(b^2 + c^2 - a^2)$, $2b^2/(c^2 + a^2 - b^2)$, $-$ 1).

The polar of a point (d, e, f) with respect to Γ has the same form as Equation (2.2) and accordingly the equation of the polar of Q with respect to Γ is

$$(b^2 + c^2 - a^2)x + (c^2 + a^2 - b^2)y + (a^2 + b^2 - c^2)z = 0. \qquad (3.4)$$

Next it emerges that the three lines AL, BN', CM' are concurrent at the point $P_1$ with co-ordinates $P_1(2a^2, a^2 - b^2 - c^2, a^2 - b^2 - c^2)$, and the three lines BM, CL', AN' are concurrent at $P_2$ with co-ordinates $P_2(b^2 - c^2 - a^2, 2b^2, b^2 - c^2 - a^2)$, and the three lines CN, AM', BL' are concurrent at $P_3$ with co-ordinates $P_3(c^2 - a^2 - b^2, c^2 - a^2 - b^2, 2c^2)$. Finally $P_1$, $P_2$, $P_3$ are collinear and lie on the line with Equation (3.4), which is the polar of Q with respect to Γ.

Note that we have now shown that triangle ABC is in perspective with each of the five triangles LMN, L'M'N', LN'M', N'ML', M'L'N, the five vertices of perspective being G, Q, $P_1$, $P_2$, $P_3$ respectively.

## 4. Definition of the points U, V, W, U', V', W', $U_2$, $U_3$, $V_3$, $V_1$, $W_1$, $W_2$

The points U, V, W are the reflections of L, M, N respectively in the sides BC, CA, AB. The points U', V', W' are the reflections of the points L', M', N' respectively in the sides BC, CA, AB. The points $U_2$, $U_3$ are the reflections of M', N' respectively in the side BC. The points $V_3$, $V_1$ are the reflections of the points N', L' respectively in the side CA and the points $W_1$, $W_2$ are the reflections of the points L', M' in the side AB. To go through all the details of obtaining the co-ordinates of these points would be both dull and unnecessary. Instead we give expressions for the reflection of a general point with co-ordinates (d, e, f) in the sides BC, CA, AB and then quote the co-ordinates of those points that are necessary to establish the features illustrated in the figure and explained in section 1.

The treatment of perpendicularity when using areal co-ordinates is tiresome. Details of how to proceed are given in Bradley [3, 4]. Following the plan explained above we find the reflection of (d, e, f) in the side BC has co-ordinates $(-d, e + d(a^2 + b^2 - c^2)/a^2, f + d(c^2 + a^2 - b^2)/a^2)$, the reflection of (d, e, f) in the side CA has co-ordinates $(d + e(a^2 + b^2 - c^2)/b^2, -e, f + e(b^2 + c^2 - a^2)/b^2)$ and the reflection of (d, e, f) in the side AB has co-ordinates $(d + f(c^2 + a^2 - b^2)/c^2, e + f(b^2 + c^2 - a^2)/c^2, -f)$.

## 5. The five Hagge circles

Using the results of Section 4 we find the co-ordinates of U to be $(a^2, 2c^2 - a^2, 2b^2 - a^2)$, the co-ordinates of V to be $(2c^2 - b^2, b^2, 2a^2 - b^2)$ and the co-ordinates of W to be $(2b^2 - c^2, 2a^2 - c^2, c^2)$. The general equation of a circle in areal co-ordinates is

$$a^2yz + b^2zx + c^2xy + (ux + vy + wz)(x + y + z) = 0. \qquad (5.1)$$

We put the co-ordinates of U, V, W into Equation (5.1) and obtain three equations for u, v, w. These are then reinstated in Equation (5.1) and simplified to obtain the Hagge circle of G, whose equation is

$$a^2(b^2 + c^2 - a^2)x^2 + b^2(c^2 + a^2 - b^2)y^2 + c^2(a^2 + b^2 - c^2)z^2 - (a^4 + (b^2 - c^2)^2)yz - (b^4 + (c^2 - a^2)^2)zx - (c^4 + (a^2 - b^2)^2)xy = 0. \qquad (5.2)$$

It may be checked that this circle passes through the orthocentre H with co-ordinates $(1/(b^2 + c^2 - a^2), 1/(c^2 + a^2 - b^2), 1/(a^2 + b^2 - c^2))$.

We now repeat this exercise to find the co-ordinates of the points U', V', W' and hence determine the equation of the Hagge circle of Q. The co-ordinates of U' turn out to be

$(1, -(a^4 - 2a^2b^2 - (b^2 - c^2)^2)/\{a^2(c^2 + a^2 - b^2)\}, -(a^4 - 2c^2a^2 - (b^2 - c^2)^2)/\{a^2(a^2 + b^2 - c^2)\})$. The co-ordinates of V' and W' may be written down by cyclic change. Using the method described above we find the Hagge circle of Q to be

$$(b^2 + c^2 - a^2)^2x^2 + (c^2 + a^2 - b^2)^2y^2 + (a^2 + b^2 - c^2)^2z^2 + (a^4 - a^2(b^2 + c^2) + 2(b^2 - c^2)^2)yz + (b^4 - b^2(c^2 + a^2) + 2(c^2 - a^2)^2)zx + (c^4 - c^2(a^2 + b^2) + 2(a^2 - b^2)^2)xy = 0. \quad (5.3)$$

Again this circle contains H. The three other Hagge circles are $UV_3W_2$, $VW_1U_3$ and $WU_2V_1$ and derive from the points $P_1$, $P_2$, $P_3$ respectively. We do not give the equations of these circles.

## 6. Three lines and three more circles

The reflections of the point L', in the sides BC, CA, AB, are the points U', $V_1$, $W_1$ respectively. These are known to lie on a line through the orthocentre H, sometimes known as the Double – Simson line of L'. Similarly the points V', $W_2$, $U_2$, H are collinear as are W', $U_3$, $V_3$, H.

The other three circles in the figure (besides Γ) are the circles BHC, CHA and AHB. These are well-known circles. The equation of CHA, for example, is

$$(c^2 + a^2 - b^2)y^2 + (c^2 - b^2)yz - b^2zx + (a^2 - b^2)xy = 0. \quad (6.1)$$

The interesting result is that this circle passes through all the points V, V', $V_1$ and $V_3$. The co-ordinates of V and V' have already been given, those of $V_1$ are (x, y, z), where

$$x = \frac{a^2 + 3(b^2 - c^2)}{c^2 + a^2 - b^2}, y = -\frac{2b^2}{c^2 + a^2 - b^2},$$
$$z = -2(a^4 - 3a^2c^2 - b^2(b^2 - c^2))/\{(a^2 + b^2 - c^2)(c^2 + a^2 - b^2)\} \quad (6.2)$$

and the co-ordinates of $V_3$ are (x, y, z), where

$$x = \frac{2\left(a^2(b^2 - 3c^2) - (b^2 + c^2)(b^2 - c^2)\right)}{\{(a^2 - b^2 - c^2)(c^2 + a^2 - b^2)\}}, y = -\frac{2b^2}{c^2 + a^2 - b^2}, z = \frac{(3b^2 + c^2 - 3a^2)}{(c^2 + a^2 - b^2)}. \quad (6.3)$$

Similarly circle BHC contains U, U', $U_2$, $U_3$ and circle AHB contains the points W, W', $W_1$, $W_2$.

## 7. When the starting point is H or K

When the starting point is the orthocentre H, then the points U, V, W coincide with H and one of the Hagge circles becomes a point. Circles such as $UV_2W_3$ are ill-defined, though, of course, the circle $HV_2W_3$ can be drawn.

When the starting point is the symmedian point K, then, as D, A, L are collinear, L, L' coincide, as do M, M' and N, N'. There are only four Hagge circles as Q coincides with K. It is interesting however, that triangle ABC is in perspective with the four triangles LMN, LNM, NML and MLN. The latter three provide what is called a triple reverse perspective and triangles ABC and LMN become triangles that define the Brocard porism. The line $P_1P_2P_3$ is now the polar of K and it contains three further interesting points AB^DE, BC^EF and CA^FD.

Flat 4,
Terrill Court,
12/14 Apsley Road,
BRISTOL BS8 2SP